\documentclass[12pt]{article}
\usepackage{fullpage}
\usepackage{amsmath}

\begin{document}
\title{Asymptotic series for Hofstadter's figure-figure sequences}
\author{
Beno{\^{\i}}t Jubin\footnote{
Mathematics Research Unit,
University of Luxembourg,
6 rue Coudenhove-Kalergi,
L-1359 Luxembourg City,
Grand-Duchy of Luxembourg
---
benoit.jubin@uni.lu}}
\date{28 May, 2013}
\maketitle
\begin{abstract}
We compute asymptotic series for Hofstadter's figure-figure sequences.
\end{abstract}

\section{Introduction}
We consider partitions of the set of strictly positive integers into two subsets such that one set, $B$, consists of the differences of consecutive elements of the other set, $A$, and a given difference appears at most once. There are many such partitions. We call $a$ the (strictly increasing) sequence enumerating $A$, and $b$ the (injective) sequence of its first differences, both with offset 1. Hofstadter's figure-figure sequences are the sequences $a$ and $b$ corresponding to the partition with the set $A$ lexicographically minimal. This is equivalent to $b$ being increasing. The sequences read
\begin{align*}
a_n &= 1, 3, 7, 12, 18, 26, 35, 45, 56, 69, \ldots &\text{(OEIS~A005228),}\\
b_n &= 2, 4, 5, 6, 8, 9, 10, 11, 13, 14, \ldots &\text{(OEIS~A030124).}
\end{align*}

These sequences were introduced by D. Hofstadter in \cite[p.~73]{geb}. They appear as an example of complementary sequences in~\cite{kim}. Their asymptotic behavior does not seem to be given anywhere in the literature except for the asymptotic equivalents mentioned by M. Hasler and D. Wilson in the related OEIS entries~\cite{oeis}.

We have by definition $b_n = a_{n+1} - a_n$, therefore $a_n = 1 + \sum_{k=1}^{n-1} b_k$. Since the sequence $a$ is strictly increasing, given any $n \geq 1$, there is a unique $k \geq 1$ such that $a_k - k < n \leq a_{k+1} - (k+1)$. This defines a sequence $u$ by letting $u_n$ be this $k$. Therefore,
\begin{equation}
\label{eq:inv}
a({u_n}) - u_n < n \leq a({u_n+1}) - (u_n + 1).
\end{equation}
The sequence $u$ is non-decreasing (actually, $u_{n+1}-u_n \in \{0,1\}$) and $u_1 = 1$. It reads
\begin{align*}
u_n &= 1, 2, 2, 2, 3, 3, 3, 3, 4, 4, \ldots &\text{(OEIS~A225687).}
\end{align*}
The partition condition implies
\begin{equation*}
b_n = n + u_n.
\end{equation*}
As a consequence,
\begin{equation}
\label{eq:int}
a_n = 1 + \frac{(n-1)n}{2} + \sum_{k=1}^{n-1} u_k.
\end{equation}

\section{Bounds and asymptotic equivalents}
Since $u_n \geq 1$, we have $a_n \geq \frac12 n(n+1)$.
Therefore the left inequality of~\eqref{eq:inv} implies
$\frac12 u_n(u_n+1) - u_n \leq n-1$
or
$u_n^2 - u_n - 2(n-1) \leq 0$
so $u_n \leq \frac12 + \sqrt{\frac14+2(n-1)}$ and finally
\begin{equation*}
1 \leq u_n < \sqrt{2n} + \frac12
\end{equation*}
so $n+ 1 \leq b_n < n + \sqrt{2n} + \frac12$ and
\begin{equation*}
b_n \sim n.
\end{equation*}
The upper bound on $u$ implies in turn
$a_n < 1 + \frac12 (n-1)n + \sum_{k=1}^{n-1}  (\sqrt{2k} + \frac12)$.
Since the function $\sqrt{\phantom{x}}$ is strictly increasing, we have
$\sum_{k=1}^{n-1} \sqrt{k} < \int_1^n \sqrt{x} \,dx = \frac23 ( n^{3/2} - 1)$.
Therefore
\begin{equation*}
\frac{n^2}2 + \frac{n}2 \leq a_n < \frac{n^2}2 + \frac{2^{3/2}}3 n^{3/2} - \frac13
\end{equation*}
and in particular
\begin{equation*}
a_n \sim \frac{n^2}2.
\end{equation*}

The relation $a_n < \frac{n^2}2 + \frac{2^3}3 \left( \frac{n}2 \right)^{3/2} - \frac13$ and the right inequality of~\eqref{eq:inv} imply $n < \frac{(u_n+1)^2}2 + \frac{2^{3/2}}3 (u_n+1)^{3/2} - u_n - \frac43$, which implies $u_n \to +\infty$. Therefore $2n \leq u_n^2 + O(u_n)$, but we saw that $u_n = O(\sqrt{n})$, so $O(u_n) \subseteq o(n)$, so $u_n^2 \geq 2n + o(n)$ therefore $u_n \geq \sqrt{2n} + o(\sqrt{n})$ and combining this with the above upper bound, we obtain
\begin{equation*}
u_n \sim \sqrt{2n}
\end{equation*}
and in particular $O(u_n) = O(\sqrt{n})$.

\section{Asymptotic series}
Since $a_n \sim \frac{n^2}2$, we have $a_{n+1} - a_n = O(n)$. Now~\eqref{eq:inv} gives $a(u_n) = n + O(u_n)$. On the other hand, \eqref{eq:int} gives $a_n = \frac{n^2}2 + \sum_{k=1}^{n-1} u_k + O(n)$, therefore $\frac{u_n^2}2 + \sum_{k=1}^{u_n-1} u_k = n + O(u_n)$. Since $u_n = O(\sqrt{n})$, we can increment the upper limit of the summation index by 1, and since $O(u_n) = O(\sqrt{n})$, we obtain the main relation
\begin{equation*}
\frac{u_n^2}2 + \sum_{k=1}^{u_n} u_k = n + O(\sqrt{n}).
\end{equation*}

We are now ready to prove by induction that for all $K \geq 1$, we have the asymptotic expansion
\begin{equation}
\label{eq:expan}
\boxed{
u_n = \sum_{k = 1}^K (-1)^{k+1} \frac{2^{1+(k-1)k/2}}{\prod_{j=1}^{k-1} (2^j+1)} \left( \frac{n}2 \right)^{1/2^k}
+ o\left( n^{1/2^K} \right).
}
\end{equation}
Indeed, the case $K=1$ reduces to $u_n \sim \sqrt{2n}$, which we already proved. We also prove the case $K=2$ separately since it is slightly different from the general case. We write $u_n= \sqrt{2n} + v_n$ with $v_n = o(\sqrt{n})$. We have
\[
\frac{u_n^2}2 - n = \sqrt{2n} \, v_n + \frac{v_n^2}2.
\]
We do not know \textit{a priori} that $v_n^2 = O(\sqrt{n})$, and that is why we have to prove this case separately. We also have
\[
\sum_{k=1}^{u_n} u_k
= \sqrt{2} \, \sum_{k=1}^{u_n} \sqrt{k} + \sum_{k=1}^{u_n} v_k
= \frac{2^{3/2}}3 {u_n}^{3/2} + o\big(O(u_n)^{3/2}\big) + \sum_{k=1}^{u_n} v_k
\]
but $\sum_{k=1}^{u_n} v_k = o\big(O(\sqrt{n})^{3/2}\big)$, therefore
\[
\frac{u_n^2}2 + \sum_{k=1}^{u_n} u_k - n = \sqrt{2n} \, v_n + \frac{v_n^2}2 + \frac{2^{9/4}}3 {n}^{3/4} + o(n^{3/4})
\]
and this has to be $O(\sqrt{n})$ by the main relation, therefore, dividing by $\sqrt{2n}$, we obtain
\[
v_n + \frac{2^{7/4}}3 {n}^{1/4} = o(n^{1/4}) + o(v_n)
\]
so
$v_n \sim - \frac{2^2}3 \left( \frac{n}2\right)^{1/4}$, as wanted.

Now, suppose the expansion holds for some $K \geq 2$. We prove it for $K+1$. It will be convenient to denote the coefficients of the expansion by
\[
\alpha_k = (-1)^{k+1} \frac{2^{1+(k-1)k/2}}{\prod_{j=1}^{k-1} (2^j+1)}.
\]
We write $v_n = o\left( n^{1/2^K} \right)$ for the remainder in~\eqref{eq:expan}. Then
\[
\frac{u_n^2}2 - n
= 2 \, \sum_{k=2}^K \alpha_k  \left( \frac{n}2 \right)^{1/2+1/2^k}
+ \sqrt{2n} \, v_n
+ O(\sqrt{n})
\]
and
\begin{align*}
\sum_{k=1}^{u_n} u_k
& = \sum_{k = 1}^K 2 \frac{2^k}{2^k+1} \alpha_k \left( \frac{u_n}2 \right)^{1+1/2^k}
+ o\left( {u_n}^{1+1/2^K} \right)\\
& = \sum_{k = 1}^K \frac{2^{k+1}}{2^k+1} \alpha_k \left( \frac{n}2 \right)^{1/2+1/2^{k+1}}
+ o\left( n^{1/2+1/2^{K+1}} \right)
\end{align*}
therefore
\begin{align*}
\frac{u_n^2}2 + \sum_{k=1}^{u_n} u_k - n
& = 2 \, \sum_{k=2}^K \alpha_k  \left( \frac{n}2 \right)^{1/2+1/2^k}
+ \sum_{k = 1}^K \alpha_k \frac{2^{k+1}}{2^k+1} \left( \frac{n}2 \right)^{1/2+1/2^{k+1}}\\
& \phantom{+} + \sqrt{2n} \, v_n + o\left( n^{1/2+1/2^{K+1}} \right)\\
& = \alpha_K \frac{2^{K+1}}{2^K+1} \left( \frac{n}2 \right)^{1/2+1/2^{K+1}}
+ 2 \left( \frac{n}2 \right)^{1/2} \, v_n + o\left( n^{1/2+1/2^{K+1}} \right)
\end{align*}
since the terms in the sums cancel out except for the last in the second sum. This expression has to be $O(\sqrt{n})$ by the main relation, so $v_n \sim - \frac{2^{K}}{2^K+1}  \alpha_K \left( \frac{n}2 \right)^{1/2^{K+1}}$, as wanted.

From the expansion of $u_n$, we find that of $b_n = n + u_n$, and that of $a_n$ by term-by-term integration. We obtain
\begin{equation*}
b_n
= n
+ \sum_{k = 1}^K (-1)^{k+1} \frac{2^{1+(k-1)k/2}}{\prod_{j=1}^{k-1} (2^j+1)} \left( \frac{n}2 \right)^{1/2^k}
+ o\left( n^{1/2^K} \right)
\end{equation*}
and
\begin{equation*}
a_n
= \frac{n^2}2
+ \sum_{k = 1}^K (-1)^{k+1} \frac{2^{k(k+1)/2}}{\prod_{j=1}^{k} (2^j+1)} \left( \frac{n}2 \right)^{1+1/2^k}
+ o\left( n^{1+1/2^K} \right).
\end{equation*}

\section{Acknowledgment}
I would like to thank Neil J. A. Sloane for useful comments on a first version of this article, and Clark Kimberling and Maximilian Hasler for their advice. I thank the Luxembourgish FNR for support via the AFR Postdoc Grant Agreement PDR~2012-1.

\bigskip
\hrule
\bigskip

\noindent 2010 {\it Mathematics Subject Classification:}
Primary 41A60.

\noindent \emph{Keywords:}
Hofstadter sequence, asymptotic series.

\bigskip
\hrule
\bigskip

\noindent (Concerned with sequences A005228, A030124, and A225687.)
\end{document}